\documentclass[11pt]{article}
\usepackage{geometry}                
\usepackage{graphicx}
\usepackage{amssymb}
\usepackage{amsmath}
\usepackage{epstopdf}
\usepackage{mathtools}
\usepackage{tikz}
\usepackage{placeins}

\usepackage{mathrsfs}

\newtheorem{theorem}{Theorem}

\newtheorem{corollary}{Corollary}

\newtheorem{proposition}{Proposition}

\begin{document}

\title {Neighborhood reconstruction and cancellation of graphs}
\author{Richard H.~Hammack\\
Department of Mathematics and Applied Mathematics\\ 
       Virginia Commonwealth University\\
       Richmond, VA 23284-2014, USA\\
       {\tt rhammack@vcu.edu} \\
 \\
 Cristina Mullican\\
 Department of Mathematics\\
 Boston College\\   
 Boston, MA   02467-3806, USA\\
{\tt mullicac@bc.edu } 
}
\date{}
\maketitle

\begin{quote}
{\small {\bf Abstract.} 
We connect two seemingly unrelated problems in graph theory.

Any graph $G$ has an associated {\it neighborhood multiset} $\mathscr{N}(G)= \{N(x) \mid x\in V(G)\}$
whose elements are precisely the open vertex-neighborhoods of $G$. In general there exist non-isomorphic graphs
$G$ and $H$ for which $\mathscr{N}(G)=\mathscr{N}(H)$. The {\it neighborhood reconstruction problem}
asks the conditions under which $G$ is uniquely reconstructible from its neighborhood multiset, that is,
the conditions under which $\mathscr{N}(G)=\mathscr{N}(H)$ implies $G\cong H$. Such a graph is said to be {\it neighborhood-reconstructible}.

The {\em cancellation problem} for the direct product of graphs seeks the conditions under which
$G\times K\cong H\times K$ implies $G\cong H$. Lov\'asz proved that this is indeed the case if $K$ is not bipartite.
A second instance of the cancellation problem asks for
conditions on $G$ that assure $G\times K\cong H\times K$ implies $G\cong H$ for any bipartite graph $K$ with $E(K)\neq \emptyset$.
A graph $G$ for which this is true is called a {\it cancellation graph}.

We prove that the neighborhood-reconstructible graphs
are precisely the cancellation graphs.
We also present some new results on cancellation graphs, which have
corresponding implications for neighborhood reconstruction.
We are particularly interested in the (yet-unsolved) problem of
finding a simple structural characterization of cancellation graphs
(equivalently, neighborhood-reconstructible graphs).
 }
\end{quote}

\medskip 

\section{Preliminaries}
For us, a graph $G$ is a symmetric relation $E(G)$ on a finite
{\it vertex set} $V(G)$.
An an {\it edge} $(x,y)\in E(G)$ is denoted $xy$. A {\it loop} is
a reflexive edge $xx$. The {\it open neighborhood} of a vertex
$x\in V(G)$ is the set $N_G(x)=\{y\in V(G)\mid xy\in E(G)\}$,
which we may denote as $N(x)$ when this is unambiguous.
Notice that $x\in N_G(x)$ if and only if $xx\in E(G)$, in which
case we say {\it there is a loop at} $x$.

In this paper we are careful to distinguish between graph
equality and isomorphism. The statement $G=H$ means $V(G)=V(H)$
and $E(G)=E(H)$. By $G\cong H$ we mean that $G$ and $H$ are
isomorphic. An isomorphism from $G$ to itself is called an
{\it automorphism} of~$G$. The group of all automorphisms of
$G$ is denoted $\mbox{Aut}(G)$. An automorphism of order 2 is
called {\it involution}.
A {\it homomorphism} $G\to H$ is a map $\varphi:V(G)\to V(H)$
for which $xy\in E(G)$ implies $\varphi(x)\varphi(y)\in E(H)$.

\begin{figure}[h]
\centering
\begin{tikzpicture}[style=thick, scale=1]
\def\vr{0.08}
\draw (-5,1)--(-5,2);
\draw (-4,0.4)--(-2,0.4);
\draw (-4,2)--(-2,2)--(-3,1)--(-4,2)--(-2,1)--(-3,2)--(-4,1)--(-2,2);
\draw (-5,2)..controls +(0.7,0.7) and +(-0.7,0.7)..+(0,0);
\draw (-4,2)..controls +(0,0.8) and +(0,0.8)..+(2,0);
\draw (-4,0.4)..controls +(0,-0.7) and +(0,-0.7)..+(2,0);
\foreach \x in {-5,-4,-3,-2} \foreach \y in {0.4,1,2} \draw [fill=white] (\x,\y) circle (\vr);
\draw [white,fill=white] (-5,0.4) circle (2*\vr);
\draw (-5.3,1.5) node {$H$};
\draw (-1.5,.4) node {$G$};
\draw (-3,2.25) node {$G\times H$};

\draw (1,0.4)--(2,0.4)--(3,0.4);
\draw (4,0.4)--(5,0.4)--(6,0.4);
\draw (2,0.4)..controls +(.5,.5) and +(-.5,.5)..+(2,0);
\draw (3,0.4)..controls +(.5,-.5) and +(-.5,-.5)..+(2,0);
\draw [densely dashed] (1,1)--(2,2)--(3,1)--(5,2)--(4,1)--(2,2);
\draw [densely dashed] (5,2)--(6,1);
\draw (1,2)--(2,1)--(3,2)--(5,1)--(4,2)--(2,1);
\draw (5,1)--(6,2);
\draw (0,1)--(0,2);
\foreach \x in {1,2,3,4,5,6} \foreach \y in {0.4,1,2} \draw  [fill=white] (\x,\y) circle (\vr);
\foreach \y in {1,2} \draw [fill=white] (0,\y) circle (\vr);
\draw (0,2.5) node {$K_2$};
\draw (6.5,0.4) node {$G$};
\draw (3.5, 2.5) node {$K_2\times G$};
\end{tikzpicture}
\caption{Examples of direct products}
\label{FIG:directp}
\end{figure}
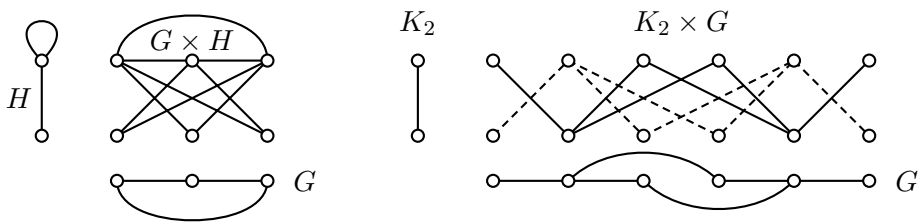

The {\it direct product} of two graphs $G$ and $H$ is the graph $G\times H$ with vertices $V(G)\times V(H)$ and edges
$E(G\times H)=\{(x,x')(y,y')\mid xy\in E(G) \mbox{ and } x'y'\in E(H)\}$.
See Figure \ref{FIG:directp}.  

We assume our reader is at least somewhat familiar with
direct products. See~\cite{HIK} for a survey.
The direct product is associative,
commutative, and distributive in the sense that
$G\times (H+K)=G\times H + G\times K$, where $+$ represents
disjoint union. Weichsel's theorem~\cite[Theorem~5.9]{HIK}  states that $G\times H$ is
connected if and only if both $G$ and $H$ are connected and at least one of them has an odd cycle. If $G$ and $H$ are both
connected and
bipartite, then $G\times H$ has exactly two components. 
In particular, if $G$ is bipartite, then $G\times K_2=G+G$,
as illustrated on the right of Figure~\ref{FIG:directp}.

\section{Neighborhood reconstruction}

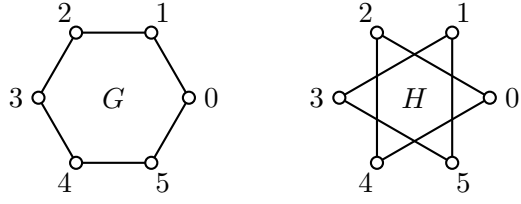
\begin{figure}[h]
\centering
\def\vr{0.08}
\begin{tikzpicture}[style=thick,scale=1]
\foreach \x in {0,1,2,3,4,5} \draw (\x*60:1)--(\x*60+60:1);
\foreach \x in {0,1,2,3,4,5} \draw (\x*60:1) [fill=white] circle (\vr);
\foreach \x in {0,1,2,3,4,5} \draw (\x*60:1.3) node {$\x$};
\draw (0,0) node {$G$};

\foreach \x in {0,1,2,3,4,5} \draw (4,0) +(\x*60:1)--+(\x*60+120:1);
\foreach \x in {0,1,2,3,4,5} \draw (4,0) +(\x*60:1) [fill=white] circle (\vr);
\foreach \x in {0,1,2,3,4,5} \draw (4,0) +(\x*60:1.3) node {$\x$};
\draw (4,0) node {$H$};
\end{tikzpicture}
\caption{Here $G\not\cong H$ but $\mathscr{N}(G)=\mathscr{N}(H)
=\{\{0,2\},\{2,4\},\{0,4\},\{1,3\},\{3,5\},\{1,5\} \}$}
\label{HEX}
\end{figure}

Any graph $G$ has an associated {\it neighborhood multiset}
$\mathscr{N}(G)= \{N_G(x) \mid x\in V(G)\}$ whose elements are
precisely the open neighborhoods of $G$. It is possible that
$G\not\cong H$ but nonetheless $\mathscr{N}(G)=\mathscr{N}(H)$,
as illustrated in Figure~\ref{HEX}.

Two types of questions have been asked about neighborhood multisets. Given a set $V$
and a multiset $\mathscr{N}=\{N_1,N_2,\ldots,N_n\}$ of subsets of $V$, we may ask
if there is a graph $G$ on $V$ for which $\mathscr{N}(G)=\mathscr{N}$.
Let us call this the {\em neighborhood realizability problem.}
Aigner and Triesch \cite{AT} attribute this problem to S\'{o}s, and show that it is NP-complete.

On the other hand, the {\it neighborhood reconstruction problem} asks whether  a given
graph $G$ can be reconstructed from the information in $\mathscr{N}(G)$, that is, whether
$\mathscr{N}(G)=\mathscr{N}(H)$ implies $G\cong H$. If this is the case we say that
$G$ is {\it neighborhood reconstructible}. Figure~\ref{HEX} shows that the hexagon is
not neighborhood reconstructible. 

Aigner and Triesch~\cite{AT} note that the problem of deciding whether a graph is
neighborhood reconstructible is NP-complete.

We now adapt their approach to describe for
given $G$ all those graphs $H$ for which $\mathscr{N}(G)=\mathscr{N}(H)$.
Given a permutation $\alpha$ of
$V(G)$ we define $G^\alpha$ to be the digraph on $V(G)$ 
with an arc directed from $x$ to $\alpha(y)$ whenever $xy\in E(G)$.
(In general, we denote an arc directed from $u$ to $v$ as an ordered list $uv$,
with the understanding that it points from the left vertex $u$ to the right vertex $v$.
Thus the arc set of $G^\alpha$ is $E(G^\alpha)=\{x\,\alpha(y)\mid xy\in E(G)\}$.)
Even though $G$ is a graph (i.e. the edge relation is symmetric), $G^\alpha$ may not be a graph. In fact, $G^\alpha$ is a graph if and only if $\alpha$ has the
property that $xy\in E(G)\Longleftrightarrow \alpha(x)\alpha^{-1}(y)\in E(G)$. Indeed, if $G^\alpha$ is a graph, then


\[\begin{array}{ccccccc}
xy\in E(G) &\Longleftrightarrow&
yx\in E(G) &\Longleftrightarrow&
y\alpha(x)\in E(G^\alpha)
&& \\
& & & \Longleftrightarrow & \alpha(x)y\in E(G^\alpha) &\Longleftrightarrow &
\alpha(x)\alpha^{-1}(y)\in E(G) .
\end{array}\]

Conversely, if $\alpha$ obeys $xy\in E(G)\Longleftrightarrow \alpha(x)\alpha^{-1}(y)\in E(G)$,
then $G^\alpha$ is a graph because
\[\begin{array}{ccccccc}
xy\in E(G^\alpha) &\Longleftrightarrow&
x \alpha^{-1}(y)\in E(G) &\Longleftrightarrow&
\alpha^{-1}(y)x\in E(G)
&& \\
& & & \Longleftrightarrow & y\alpha^{-1}(x)\in E(G) &\Longleftrightarrow &
yx\in E(G^\alpha) .
\end{array}\]

A map $\alpha$ with the above properties is called an {\it anti-automorphism}
in~\cite{HIK} and~\cite{hammack2009}.

To summarize, an {\it anti-automorphism} of a graph $G$ is a bijection
$\alpha:V(G)\to V(G)$ for which $xy\in E(G)$ if and only if
$\alpha(x)\alpha^{-1}(y)\in E(G)$. Given an anti-automorphism $\alpha$ of
$V(G)$ we have a graph $G^\alpha$ on the same vertex set as $G$, but with
\[E(G^\alpha)=\{ x \alpha(y) \mid xy\in E(G)\}.\]
Notice that this means $N_{G}(y)=N_{G^\alpha}\big(\alpha(y)\big)$,
and therefore
\begin{equation}
\mathscr{N}(G)=\mathscr{N}(G^\alpha).
\label{Eqn:N=N}
\end{equation}

For example, consider the hexagon $G$ in Figure~\ref{HEX}, and let
$\alpha$ be the antipodal map that rotates it $180^\circ$ about its center.
Then $\alpha$ is an anti-automorphsm (it also happens to be an
automorphism) and $G^\alpha=H$ is the union of two
triangles shown on the right of Figure~\ref{HEX}. 

\begin{figure}[h]
\centering
\begin{tikzpicture}[style=thick,scale=0.7]
\def\vr{0.115} 
\draw (0,0)--(2,0)--(2,2)--(0,2)--(0,0) (0,0)--(2,2);
\draw (0,2)..controls +(-1.2,0) and +(0,1.2)..+(0,0);
\draw (2,0)..controls +(1.2,0) and +(0,-1.2)..+(0,0);
\draw (1.5,1) node {$G$}; \draw (2.8,1)  node {$\alpha$};
\draw [->,>=stealth,densely dashed,shorten >=0.15cm] (2,2) arc(45:-45:1.41);
\draw [->,>=stealth,densely dashed,shorten >=0.15cm] (2,0) arc(-45:-135:1.41);
\draw [->,>=stealth,densely dashed,shorten >=0.15cm] (0,0) arc(-135:-225:1.41);
\draw [->,>=stealth,densely dashed,shorten >=0.15cm] (0,2) arc(135:45:1.41);
\foreach \x in {0,2} \foreach \y in {0,2} \draw (\x,\y) [fill=white] circle (\vr);
\end{tikzpicture}
\hspace{0.2in}
\begin{tikzpicture}[style=thick,scale=0.7]
\def\vr{0.115} 
\draw (0,0)--(2,0);
\draw (2,2)--(0,2);
\draw (0,0)--(2,2);
\draw (0,2)--(2,0);
\draw (0,0)..controls +(-1.2,0) and +(0,-1.2)..+(0,0);
\draw (0,2)..controls +(-1.2,0) and +(0,1.2)..+(0,0);
\draw (2,2)..controls +(1.2,0) and +(0,1.2)..+(0,0);
\draw (2,0)..controls +(1.2,0) and +(0,-1.2)..+(0,0);
\foreach \x in {0,2} \foreach \y in {0,2} \draw (\x,\y) [fill=white] circle (\vr);
\draw (1.7,1) node {$G^\alpha$}; 
\end{tikzpicture}
\caption{A graph $G$, an anti-automorphism $\alpha$, and the
corresponding graph $G^\alpha$. }
\label{FigSquare}
\end{figure}
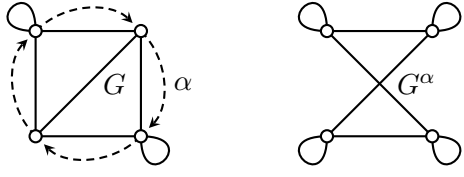

Figure~\ref{FigSquare} shows a second example. The $90^\circ$
rotation $\alpha$ of $V(G)$ is an anti-automorphism, and $G^\alpha$
is shown on the right.  Notice that $\mathscr{N}(G)=\mathscr{N}(G^\alpha)$.
These examples illustrate our
first proposition, which in essence was noted in~\cite{AT}, in the context of loopless
graphs.\footnote{The article~\cite{AT} differs slightly from our current setting. What we here call
an {\it anti-automorphism} plays the role of an {\it admissible map}
in~\cite{AT}. Admissible maps coincide with our anti-automorphisms,
except that they have an additional condition that 
assures $G^\alpha$ is loopless. Thus the definition of an
anti-automorphism is weaker than that of an admissible map.}

\begin{proposition}
If $G$ and $H$ are two graphs on the same vertex set, then 
$\mathscr{N}(G)=\mathscr{N}(H)$ if and only if $H=G^\alpha$ for some
anti-automorphism of $G$.
\label{PROP:neighborhood}
\end{proposition}

\noindent
Proof: 
If $H=G^\alpha$ for some anti-automorphism $\alpha$ of $G$,
then $\mathscr{N}(G)=\mathscr{N}(H)$ by Equation~\ref{Eqn:N=N}.

Conversely, let
$G$ and $H$ have vertex set $V$, and suppose
$\mathscr{N}(G)=\mathscr{N}(H)$.
Then there is a permutation $\alpha$ of $V$ with
$N_G(x)=N_H\big(\alpha(x)\big)$ for all $x\in V$, so also
$N_H(x)=N_G\big(\alpha^{-1}(x)\big)$.
Note that $\alpha$ is an anti-automorphism of $G$ because
\begin{align*}
xy\in E(G)&\Longleftrightarrow y\in N_G(x)\\
&\Longleftrightarrow y\in N_H\big(\alpha(x)\big)\\
&\Longleftrightarrow y\alpha(x)\in E(H)\\
&\Longleftrightarrow \alpha(x)\in N_H(y)\\
&\Longleftrightarrow \alpha(x)\in N_G\big(\alpha^{-1}(y)\big)\\
&\Longleftrightarrow \alpha(x)\alpha^{-1}(y)\in E(G).
\end{align*}

To verify $H=G^\alpha$, observe that  
\begin{align*} 
xy\in E(H)&\Longleftrightarrow x\in N_H(y)\\
&\Longleftrightarrow x\in N_G\big(\alpha^{-1}(y)\big)\\
&\Longleftrightarrow x\alpha^{-1}(y)\in E(G)\\
&\Longleftrightarrow x\alpha(\alpha^{-1}(y))=xy\in E(G^\alpha).
\tag*{\rule{7pt}{7pt}}
\end{align*}


\medskip
Let's pause to elaborate on the notion of neighborhood
reconstructibility. As noted above, $G$ is neighborhood-reconstructible if for any $H$ with
$\mathscr{N}(G)=\mathscr{N}(H)$ it necessarily follows that
$G\cong H$. Define $G$ to be {\it strongly neighborhood-reconstructible} if $\mathscr{N}(G)=\mathscr{N}(H)$ implies
$G=H$.

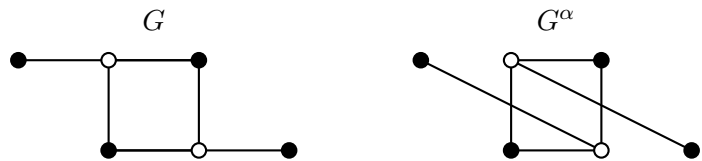
\begin{figure}[t]
\centering
\begin{tikzpicture}[style=thick,scale=0.6]
\def\vr{0.16}
\draw (3,-1)--(-1,-1);
\draw (1,1)--(-3,1);
\draw (-1,-1)--(-1,1)--(1,1)--(1,-1)--(-1,-1);
\draw (3,-1) [fill=black] circle (\vr);
\draw (1,1) [fill=black] circle (\vr);
\draw (-1,-1) [fill=black] circle (\vr);
\draw (-3,1) [fill=black] circle (\vr);
\draw (-1,1) [fill=white] circle (\vr);
\draw (1,-1) [fill=white] circle (\vr);
\draw (0,1.9) node {$G$};
\end{tikzpicture}
\hspace{0.5in}
\begin{tikzpicture}[style=thick,scale=0.6]
\def\vr{0.16}
\draw (3,-1)--(-1,1);
\draw (1,-1)--(-3,1);
\draw (-1,-1)--(-1,1)--(1,1)--(1,-1)--(-1,-1);
\draw (3,-1) [fill=black] circle (\vr);
\draw (1,1) [fill=black] circle (\vr);
\draw (-1,-1) [fill=black] circle (\vr);
\draw (-3,1) [fill=black] circle (\vr);
\draw (-1,1) [fill=white] circle (\vr);
\draw (1,-1) [fill=white] circle (\vr);
\draw (0,1.9) node {$G^\alpha$};
\end{tikzpicture}
\caption{Rotation $\alpha$ of $G$ by $180^\circ$ is an involution, and hence also an anti-automorphism.
Notice that $\mathscr{N}(G)=\mathscr{N}(G^\alpha)$.
Here $G\ne G^\alpha$, but $G\cong G^\alpha$.}
\label{FigReconstruct}
\end{figure}

The example in Figure~\ref{FigReconstruct} should clarify
the distinction.
Clearly $G\ne G^\alpha$, though $G\cong G^\alpha$ and $\mathscr{N}(G)=\mathscr{N}(G^\alpha)$.
In fact, the results of Section~\ref{Section:Further} will show
that this graph $G$ is neighborhood-reconstructible. Hence it
is neighborhood-reconstructible but not strongly
neighborhood-reconstructible.

For a simple (but not completely trivial) example of a graph that is strongly
neighborhood-reconstructible, let $G$ be an edge $ab$ with
a loop at $a$. It is straightforward that $G$ 
with $E(G)=\{ab, aa\}$ is the only
graph that can be reconstructed from $\mathscr{N}(G)$.

We close this section with an immediate corollary of
Proposition~\ref{PROP:neighborhood}.

\begin{corollary}
A graph $G$ is neighborhood-reconstructible if and only if
$G\cong G^\alpha$ for every anti-automorphism $\alpha$ of $G$.
\label{Cor:NbhdAnti}
\end{corollary}

Observe that Proposition~\ref{PROP:neighborhood} also implies that $G$
is strongly neighborhood-reconstructible if and only if
$G= G^\alpha$ for every anti-automorphism $\alpha$ of $G$.
We can further refine this by forming an equivalence relation
$R$ on $V(G)$ by declaring $xRy$ if and only if $N(x)=N(y)$.
The proof of the next corollary is straightforward from definitions.

\begin{corollary}
A graph $G$ is strongly neighborhood-reconstructible if and only if
its anti-automorphisms  are precisely the permutations
of $V(G)$ that preserve the $R$-equivalence classes
of $V(G)$. (That is, $\alpha(X)=X$ for each $R$-equivalence class $X$.)
\end{corollary}

Although Corollary~\ref{Cor:NbhdAnti} characterizes neighborhood-reconstructible
graphs, we certainly cannot regard it as a simple
characterization, as finding all anti-automorphisms of $G$
promises to be quite difficult in general, let alone
deciding if $G\cong G^\alpha$ for all of them. However,
it does provide a link to cancellation, which we now explore.

\section{Cancellation}

Lov\'{a}sz \cite[Theorem 9]{lovasz} proved that if a graph $K$ has an odd cycle, then
$G\times K\cong H\times K$ implies $G\cong H$. In such a situation
we say that {\it cancellation holds}.

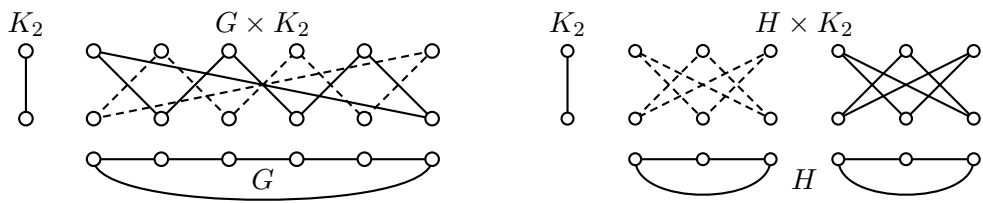
\begin{figure}[h]
\centering
\begin{tikzpicture}[style=thick,scale=0.9]
\def\vr{0.088}
\draw (0,1)--(0,2);
\draw (1,0.4)--(6,0.4);
\draw [densely dashed](1,1)--(2,2)--(3,1)--(4,2)--(5,1)--(6,2)--(1,1);
\draw (1,2)--(6,1)--(5,2)--(4,1)--(3,2)--(2,1)--(1,2);
\draw (1,0.4)..controls +(0,-0.8) and +(0,-0.8)..+(5,0);
\foreach \x in {1,2,3,4,5,6} \foreach \y in {0.4,1,2} \draw [fill=white] (\x,\y) circle (0.1);
\draw [fill=white] (0,1) circle (0.1);
\draw [fill=white] (0,2) circle (0.1);
\draw (0,2.4) node {$K_2$};
\draw (3.5,0.1) node {$G$};
\draw (3.5,2.4) node {$G\times K_2$};

\draw (8,1)--(8,2);
\draw (9,0.4)--(11,0.4);
\draw (12,0.4)--(14,0.4);
\draw [densely dashed](9,1)--(10,2)--(11,1)--(9,2)--(10,1)--(11,2)--(9,1);
\draw (12,1)--(13,2)--(14,1)--(12,2)--(13,1)--(14,2)--(12,1);
\draw (9,0.4)..controls +(0,-0.7) and +(0,-0.7)..+(2,0);
\draw (12,0.4)..controls +(0,-0.7) and +(0,-0.7)..+(2,0);
\foreach \x in {9,10,11,12,13,14} \foreach \y in {0.4,1,2} \draw [fill=white] (\x,\y) circle (\vr);
\draw [fill=white] (8,1) circle (\vr);
\draw [fill=white] (8,2) circle (\vr);
\draw (8,2.4) node {$K_2$};
\draw (11.5,0.1) node {$H$};
\draw (11.5,2.4) node {$H\times K_2$};
\draw (15,0) [white] circle (\vr);
\end{tikzpicture}
\caption{Failure of cancellation: $G\times K_2\cong H\times K_2$
but $G\not\cong H$.}
\label{dpexample}
\end{figure}

Cancellation may fail if
$K$ is bipartite. For example, consider graphs $G$ and $H$ from Figure~\ref{HEX}.
In Figure~\ref{dpexample} we see that
$G\times K_2\cong H\times K_2$, as both products are isomorphic to
two copies of a hexagon, but cancellation fails because
$G\not\cong H$. Recall that $H=G^\alpha$ where $\alpha$ is
the antipodal map of $G$. Thus we have 
$G\times K_2\cong G^\alpha\times K_2$.

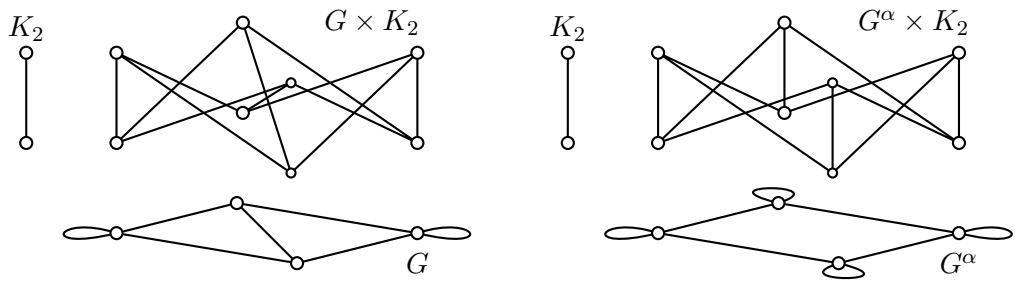
\begin{figure}[h]
\centering
\begin{tikzpicture}[style=thick,scale=0.4]
\def\vr{0.2} 
\draw (0,0)--+(6,-1)--+(10,0)--+(4,1)--+(0,0) +(4,1)--+(6,-1);
\draw (0,0)  +(-3,3)--+(-3,6);
\draw (0,0) +(0,0)..controls +(-2,0.66) and +(-2.66,-0.66)..+(0,0);
\draw (10,0)..controls +(2,0.66) and +(2.66,-0.66)..+(0,0);
\draw (0,0) +(0,6)--+(5.8,2)--+(10,6)--+(4.2,4)--+(0,6);
\draw (0,0) +(0,3)--+(5.8,5)--+(10,3)--+(4.2,7)--+(0,3);
\draw (0,0) +(0,3)--+(0,6);
\draw (0,0) +(10,3)--+(10,6);
\draw (0,0) +(4.2,7)--+(5.8,2);
\draw (0,0) +(4.2,4)--+(5.8,5);
\draw (10,-1) node {$G$};
\draw (-3,6.8) node {$K_2$};
\draw (8.5,7) node {$G\times K_2$};

\draw (18,0)--+(6,-1)--+(10,0)--+(4,1)--+(0,0);
\draw (18,0)  +(-3,3)--+(-3,6);
\draw (18,0) +(0,0)..controls +(-2,0.66) and +(-2.66,-0.66)..+(0,0);
\draw (28,0)..controls +(2,0.66) and +(2.66,-0.66)..+(0,0);
\draw (22,1)..controls +(2,0.66) and +(-2.66,0.66)..+(0,0);
\draw (24,-1)..controls +(-2,-0.66) and +(2.66,-0.66)..+(0,0);
\draw  (18,0) +(0,6)--+(5.8,2)--+(10,6)--+(4.2,4)--+(0,6);
\draw  (18,0) +(0,3)--+(5.8,5)--+(10,3)--+(4.2,7)--+(0,3);
\draw (18,0) +(0,3)--+(0,6);
\draw (18,0) +(10,3)--+(10,6);
\draw (18,0) +(4.2,7)--+(4.2,4);
\draw (18,0) +(5.8,2)--+(5.8,5);
\draw (28,-1) node {$G^\alpha$};
\draw (15,6.8) node {$K_2$};
\draw (26.5,7) node {$G^\alpha\times K_2$};

\foreach \x in {-3, 0,10,15,18,28} \foreach \y in {3,6} \draw (\x,\y) [fill=white] circle (\vr);
\foreach \x in { 0,10,18,28} \foreach \y in {0} \draw (\x,\y) [fill=white] circle (\vr);
\foreach \x in { 4.2, 22.2} \foreach \y in {4,7} \draw (\x,\y) [fill=white] circle (\vr);
\foreach \x in { 5.8, 23.8} \foreach \y in {2,5} \draw (\x,\y) [fill=white] circle (0.15);
\foreach \x in {4,22} \foreach \y in {1} \draw (\x,\y) [fill=white] circle (\vr);
\foreach \x in {6,24} \foreach \y in {-1} \draw (\x,\y) [fill=white] circle (\vr);
\end{tikzpicture}
\caption{Failure of cancellation: $G\times K_2\cong G^\alpha\times K_2$
but $G\not\cong G^\alpha$.}
\label{Fig:SquareK}
\end{figure}

For another example, take the graphs $G$ and $G^\alpha$ from
Figure~\ref{FigSquare}. Again, Figure~\ref{Fig:SquareK} reveals
that $G\times K_2\cong G^\alpha\times K_2$. (Each products is
isomorphic to the three-dimensional cube.)

These examples are instances of our next proposition,
which was proved in~\cite{hammack2009} and also in~\cite{thesis}.
For completeness we include an abbreviated proof.

\begin{proposition} Suppose a bipartite graph $K$  has at least one edge.
Then $G\times K\cong H\times K$ if and only if
$H\cong G^\alpha$ for some anti-automorphism $\alpha$ of $G$.
\label{PROP:permutediso}
\end{proposition}

\noindent
Proof. We use a result by
Lov\'{a}sz~\cite[Theorem 6]{lovasz}:
If there is a graph homomorphism $K'\to K$, then
$G\times K\cong H\times K$ implies $G\times K'\cong H\times K'$.

Let $G\times K\cong H\times K$.
As there is a homomorphism
$K_2\to K$, we have
$G\times K_2\cong H\times K_2$.
Take an isomorphism $\varphi:G\times K_2\to H\times K_2$.
We easily check that we
may assume $\varphi$ has form
$\varphi(g,k)=(\beta(g,k),k)$. 
(This is also a special instance of~\cite[Theorem 7]{lovasz}.)
Put $V(K_2)=\{0,1\}$.

Define bijections $\mu,\lambda:G\rightarrow H$ as $\mu(g)=\beta(g,0)$ and
$\lambda(g)=\beta(g,1)$.  First we will show $\mu^{-1}\lambda$
is an anti-automorhpism of $G$.  Then we will show $G^{\mu^{-1}\lambda}\cong H$.
Observe that
\begin{eqnarray*}
xy\in E(G) &\Longleftrightarrow & (x,1)(y,0)\in E(G\times K_2)\\
		&\Longleftrightarrow &\varphi(x,1)\varphi(y,0)\in E(H\times K_2)\\
		&\Longleftrightarrow &(\lambda(x),1)(\mu(y),0)\in E(H\times K_2)\\
		&\Longleftrightarrow &\lambda(x)\mu(y)\in E(H).
\end{eqnarray*}
A similar argument gives $xy\in E(G)\Longleftrightarrow \mu(x)\lambda(y)\in E(H)$. It follows that	
\begin{eqnarray*}
xy\in E(G) &\Longleftrightarrow&  \lambda(x)\mu(y)\in E(H)\\
		&\Longleftrightarrow& \mu^{-1}(\lambda(x))\lambda^{-1}(\mu(y))\in E(G)\\
		&\Longleftrightarrow &\mu^{-1}\lambda(x)(\mu^{-1}\lambda)^{-1}(y)\in E(G).
\end{eqnarray*}
Therefore $\mu^{-1}\lambda$ is an anti-automorphism of $G$.

Now, $\mu:G^{\mu^{-1}\lambda}\rightarrow H$ is an isomorphism because
\begin{eqnarray*}
xy\in E(G^{\mu^{-1}\lambda})	 &\Longleftrightarrow & x(\mu^{-1}\lambda)^{-1}(y)\in E(G)\\
						&\Longleftrightarrow  &  x\lambda^{-1}(\mu(y))\in E(G)\\
						&\Longleftrightarrow & \mu(x)\lambda(\lambda^{-1}\mu(y))\in E(H)\\
        &\Longleftrightarrow    &    \mu(x)\mu(y)\in E(H).
\end{eqnarray*}

Conversely, suppose $H\cong G^\alpha$ for an anti-automorphism 
$\alpha$ of $G$. Say $K$ has partite sets $X_0$ and $X_1$.
It suffices to show that $G\times K\cong G^\alpha\times K$.
Define $\Theta:G\times K\rightarrow G^\alpha\times K$ to be $$\begin{displaystyle}{\Theta(g,k)=\left\{\begin{array}{rl}
(g,k)&\mbox{if }k\in X_0\\
(\alpha(g),k)&\mbox{if }k\in X_1.\\
\end{array}\right.}\end{displaystyle}$$
To prove $\Theta$ is an isomorphism take
$(g,k)(g',k')\in E(G\times K)$. Then $k$ and $k'$ must be in
different partite sets. Say $k\in X_0$ and $k'\in X_1$. Then
\begin{align*}
(g,k)(g',k')\in E(G\times K)&\Longleftrightarrow gg'\in E(G)\text{ and }kk'\in E(K)&& \\ 
&\Longleftrightarrow g\alpha(g')\in E(G^\alpha)\text{ and }kk'\in E(K)&&(\text{def. of } G^\alpha)\\
&\Longleftrightarrow (g,k)(\alpha(g'),k')\in E(G^\alpha \times K)&&\\ 
&\Longleftrightarrow \Theta(g,k)\Theta(g',k')\in E(G^\alpha\times K).&&
\end{align*}

Similarly, if $k\in X_1$ and $k'\in X_0$ we get $(g,k)(g',k')\in E(G\times K)$ if and only if
$\Theta(g,k)\Theta(g',k')\in E(G^\alpha\times K)$.
Thus $\Theta$ is an isomorphism. \hfill \rule{7pt}{7pt}

\medskip
We define a graph $G$ to be a \textsl{cancellation graph} if
$G\times K\cong H\times K$ implies $G\cong H$ for all graphs $K$
that have at least one edge.
(We require at least one edge because if $K$ is edgeless, then
$G\times K\cong H\times K$ whenever $|V(G)|=|V(H)|$, as both
products are also edgeless.)

Note that $G^{\mbox{\scriptsize id}}=G$.
By Lov\'{a}sz's result that $G\times K\cong H\times K$ implies $G\cong H$ when $K$ has an odd cycle, we see that $G\times K\cong H\times K$ if and only if $H\cong G^{\mbox{\scriptsize id}}$
when $K$ is not bipartite.
Combining this with our
Proposition~\ref{PROP:permutediso}, we get a corollary.

\begin{corollary} A graph $G$ is a cancellation graph if and only if $G\cong G^\alpha$ for all its anti-automorphisms~$\alpha$. 
\label{COR:can}
\end{corollary}

Combining this with Corollary~\ref{Cor:NbhdAnti} yields a theorem.

\begin{theorem} A graph is a cancellation graph if and only if it is neighborhood-reconstructible.
\label{cancellationneighborhood}
\end{theorem}

\section{Further Results}
\label{Section:Further}

This section seeks structural characterizations of cancellation (hence neighborhood
reconstructible) graphs. We develop a sufficient condition for arbitrary graphs,
and a characterization for the bipartite case.

Let $\mbox{Ant}(G)$ denote the set of all anti-automorphisms of
$G$. By Corollary~\ref{Cor:NbhdAnti}, Proposition~\ref{PROP:permutediso} and Theorem~\ref{cancellationneighborhood}, a graph is 
neighborhood-reconstructible and a cancellation graph if and
only if $G\cong G^\alpha$ for each $\alpha\in\mbox{Ant}(G)$.
Therefore it is beneficial to determine the conditions
under which $G\cong G^\alpha$, and, more generally,
when $G^\alpha\cong G^\beta$.

To this end we adopt a construction from~\cite{TF2}, and apply it to our current
setting. Define the following set of pairs
of permutations of $V(G)$. Let
\[\mbox{\small
$\mbox{Aut}^{\scriptsize {\rm TF}}(G)=\big\{(\lambda,\mu) \mid
\lambda, \mu \mbox{ are permutations of $V(G)$ with
$xy\in E(G)\Longleftrightarrow \lambda(x)\mu(y)\in E(G)$}
\big\}$}.\]
Elements of $\mbox{Aut}^{\scriptsize {\rm TF}}(G)$ are called {\it two-fold automorphisms of $G$} in~\cite{TF1} and~\cite{TF2}.
Notice that $\mbox{Aut}^{\scriptsize {\rm TF}}(G)$ is non-empty because
it contains $(\mbox{id},\mbox{id})$. It is also a group under
pairwise composition, and with
$(\lambda,\mu)^{-1}=(\lambda^{-1},\mu^{-1})$.
Observe that $(\lambda,\mu)\in \mbox{Aut}^{\scriptsize {\rm TF}}(G)$ if and only if
$\lambda\big(N_G(x)\big)=N_G\big(\mu(x)\big)$ for all $x\in V(G)$.
Also $\alpha\in\mbox{Ant}(G)$ if and only if 
$(\alpha,\alpha^{-1})\in\mbox{Aut}^{\scriptsize {\rm TF}}(G)$, and 
$\alpha\in\mbox{Aut}(G)$ if and only if 
$(\alpha,\alpha)\in\mbox{Aut}^{\scriptsize {\rm TF}}(G)$.

We can think of $\mbox{Aut}^{\scriptsize {\rm TF}}(G)$ as follows:
Suppose $\lambda:V(G)\to V(G)$ is a bijection that sends
neighborhoods of $G$ to neighborhoods of $G$, that is it
``permutes'' the elements of $\mathscr{N}(G)$. Then there
must be at least one bijection $\mu:V(G)\to V(G)$ for which
$\lambda\big(N_G(x)\big)=N_G\big(\mu(x)\big)$,
and then $(\lambda,\mu)\in \mbox{Aut}^{\scriptsize {\rm TF}}(G)$.
If no two vertices of $G$ have the same neighborhood, then
there is a unique $\mu$ paired with any such $\lambda$, otherwise
there will be more than one $\mu$.
\footnote{
We remark in passing that $\mbox{Aut}^{\scriptsize {\rm TF}}(G)$ is similar
to the so-called {\it factorial} $G!$ of a graph (or digraph) $G$, as defined
in~\cite{hammackcancel} and~\cite{HIK}. The vertex set of $G!$ is the set of permutations
of $V(G)$, with an edge joining permutations $\lambda$ and $\mu$
provided $xy\in E(G)$ implies $\lambda(x)\mu(y)\in E(G)$. Thus the edge set
of $G!$ can be identified with $\mbox{Aut}^{\scriptsize {\rm TF}}(G)$.
The factorial is used in~\cite{hammackcancel} to settle the general cancellation
problem for digraphs. However, our present purposes do not require the graph
structure of the factorial, so we will phrase the discussion in terms of two-fold
automorphims.}

The group $\mbox{Aut}^{\scriptsize {\rm TF}}(G)$ acts on the set $\mbox{Ant}(G)$
as $(\lambda, \mu)\cdot\alpha = \lambda\alpha\mu^{-1}$.

\begin{proposition} Suppose $\alpha,\beta\in{\rm Ant}(G)$. 
Then $G^\alpha\cong G^\beta$ if and only if
$\alpha$ and $\beta$ are in the same ${\rm Aut}^{\scriptsize {\rm TF}}(G)$-orbit.
In particular, $G$ is neighborhood reconstructible and a cancellation graph if and only if the 
${\rm Aut}^{\scriptsize {\rm TF}}(G)$ action on ${\rm Ant}(G)$ is transitive.
\label{PROP:simeqiso}
\end{proposition}
\noindent
Proof:
Suppose $\gamma:G^\alpha\rightarrow G^\beta$ is an isomorphism.  Then
\begin{align*}
xy\in E(G) &\Longleftrightarrow x\alpha(y)\in E(G^\alpha)\\
							&\Longleftrightarrow \gamma(x)\gamma(\alpha(y))\in E(G^\beta)\\
							&\Longleftrightarrow \gamma(x)\beta^{-1}(\gamma(\alpha(y)))\in E(G).
\end{align*}
Thus $(\gamma,\,\beta^{-1}\gamma\alpha)\in  \mbox{Aut}^{\scriptsize {\rm TF}}(G)$. Also
$\gamma\alpha(\beta^{-1}\gamma\alpha)^{-1}=\beta$,
so $\alpha$ and $\beta$ are in the same $\mbox{Aut}^{\scriptsize {\rm TF}}(G)$-orbit.

Conversely, let $\alpha$ and $\beta$ be in the same  $\mbox{Aut}^{\scriptsize {\rm TF}}(G)$-orbit. Take $(\lambda,\mu)\in  \mbox{Aut}^{\scriptsize {\rm TF}}(G)$ with
$\beta=(\lambda,\mu)\cdot \alpha=\lambda\alpha\mu^{-1}$, 
so $\alpha^{-1}\lambda^{-1}=\mu^{-1}\beta^{-1}$.
Then $\lambda^{-1}:G^\beta\rightarrow G^\alpha$ is an isomorphism
because
\begin{align}
\nonumber
xy\in E(G^\beta) &\Longleftrightarrow x\beta^{-1}(y)\in E(G)\\
										&\Longleftrightarrow \lambda^{-1}(x)\mu^{-1}(\beta^{-1}(y))\in E(G)&& \tag{\mbox{\small because  $(\lambda^{-1},\mu^{-1})\in  \mbox{Aut}^{\scriptsize {\rm TF}}(G)$}}\\                                       
										&\Longleftrightarrow \lambda^{-1}(x)\alpha^{-1}(\lambda^{-1}(y))\in E(G)&& \tag{\mbox{\small because $\alpha^{-1}\lambda^{-1}=\mu^{-1}\beta^{-1}$}}\\
										&\Longleftrightarrow \lambda^{-1}(x)\lambda^{-1}(y)\in E(G^\alpha). & & \tag*{\rule{7pt}{7pt}}
\end{align}

\medskip
If $\alpha\in{\rm Ant}(G)$, then it is immediate that also $\alpha^k\in{\rm Ant}(G)$ for all integers $k$.
Moreover, because $(\alpha,\alpha^{-1})\in \mbox{Aut}^{\scriptsize {\rm TF}}(G)$ and $(\alpha,\alpha^{-1})\cdot\alpha=\alpha^{3}$, Proposition~\ref{PROP:simeqiso} yields
$G^{\alpha}\cong G^{\alpha^{3}}$. Iterating, we get a proposition.

\begin{proposition} If $\alpha\in {\rm Ant}(G)$, then $G^{\alpha}\cong G^{\alpha^{1+2n}}$ for all integers $n$. 
In particular, if $\alpha$ has odd order, then $G^\alpha\cong G$. 
\label{PROP:simplus2}
\end{proposition}


Now, if $\alpha$ has {\it even} order, we may write its order as $2^m(1+2n)$ for integers $m$ and $n$.
Then $\alpha^{1+2n}$ has order $2^m$, and by Proposition~\ref{PROP:simplus2}, $G^\alpha\cong G^{\alpha^{1+2n}}$.
Consequently we can get all $G^\alpha$, up to isomorphism, with only those anti-automorphisms whose order
is a power of~$2$. Of course this is little help in enumerating all $G^\alpha$, but it does lead to a quick
sufficient condition for a graph to be neighborhood-reconstructible.







\begin{corollary}
If a graph has no involutions, then it is  neighborhood-reconstructible, and thus also a cancellation graph.
\label{COR:can2}
\end{corollary}
\noindent
Proof:  Suppose that $G$ is not neighborhood-reconstructible. Then there is some  $\alpha\in{\rm Ant}(G)$
with $G^\alpha\not\cong G$. Proposition~\ref{PROP:simplus2} says the order $n$ of $\alpha$ is even, so
$\alpha^{n/2}$ is an involution of $G$. \hfill \rule{7pt}{7pt}



\medskip
If $G$ is bipartite, this corollary tightens to a characterization.
As a preliminary to this we claim that any anti-automorphism $\alpha$ of a bipartite graph
carries any partite set of a connected component of $G$ bijectively to a partite set
of a component of $G$. Indeed, suppose $x_0$ and $x_0'$ both belong to the same
partite set of a connected component of $G$.
Then $G$ has an
even-length path $x_0,v_1,...v_{2n+1},x_0'$. Thus the path
$\alpha(x_0),\alpha^{-1}(v_1),...,\alpha^{-1}(v_{2n+1}),\alpha(x_0')$ has even length, so
$\alpha(x_0)$ and $\alpha(x_0')$ are in the same partite set of  some component of $G$.


\medskip
\begin{proposition} A  bipartite graph is a cancellation graph (is neighborhood-reconstructible) if and only if 
it has no involution that reverses the bipartition of one of its components.
\label{PROP:biprevinv}
\end{proposition}

\noindent
Proof: Let $G$ be bipartite. Suppose $G$ has an involution $\alpha$ that
reverses the partite sets  of one of its components. Call that component $H$, and its partite sets
$X$ and $Y$. Select $x\in X$. Then $\alpha(x)\in Y$, and $H$ has an odd path
$x,x_1,x_2,x_3,\ldots, x_{2k-1},x_{2k},\alpha(x)$. Thus $G^\alpha$ has an odd walk
$x,\alpha(x_1),x_2,\alpha(x_3),\ldots, \alpha(x_{2k-1}),x_{2k},\alpha^2(x)$.
But this odd walk begins and ends at $x$, so $G^\alpha$ is not bipartite.
Consequently $G\not\cong G^\alpha$ so $G$ is neither a cancellation graph
nor neighborhood-reconstructible, by Corollaries~\ref{Cor:NbhdAnti} and~\ref{COR:can}.


Conversely, suppose $G$ has no involutions that reverse the bipartition of a component.
Say $G$ has $c$ components $H_i$, each with partite sets $V(H_i)=X_i\cup Y_i$, where
$1\leq i\leq c$. Now, we noted above that any $\alpha\in \mbox{Ant}(G)$ permutes the
set $\{X_1,Y_1,X_2,Y_2\ldots X_c,Y_c\}$. Notice that the $\alpha$-orbit of a particular $X_i$ cannot meet
the $\alpha$-orbit of the corresponding~$Y_i$. The reason is that we'd then have $\alpha^k(X_i)=Y_i$ for some power $k$. From this we could concoct an involution $\sigma$ of $G$ that reverses the bipartition of
$H_i$ by simply declaring
$$\begin{displaystyle}{\sigma(x)=\left\{\begin{array}{cl}
x&\mbox{if }x\in V(G)-V(H_i)\\
\alpha^{k}(x)&\mbox{if }x\in X_i\\
\alpha^{-k}(x)&\mbox{if }x\in Y_i.\\
\end{array}\right.}\end{displaystyle}$$
Since no such involution exists, the $\alpha$-orbit of a $X_i$ never meets
the $\alpha$-orbit of $Y_i$.

Therefore, given any $\alpha\in \mbox{Ant}(G)$, we may assume (by interchanging
the labels $X_i$ and $Y_i$ as appropriate) that $\alpha$ sends each $X_i$ to some
$X_j$, and it sends each $Y_k$ to some $Y_\ell$.
Define a bipartition
$V(G)=X\cup Y$, where $X=\bigcup X_i$, and $Y=\bigcup Y_i$.
By construction we have $\alpha(X)=X$ and $\alpha(Y)=Y$.
Now form a map $\mu:G\to G^\alpha$ as
$$\begin{displaystyle}{\mu(x)=\left\{\begin{array}{ll}
\alpha(x)&\mbox{if }x\in X\\
x &\mbox{if }x\in Y.\\
\end{array}\right.}\end{displaystyle}$$
That this is an isomorphism follows immediately from the definition of $G^\alpha$ and
the anti-automorphism property of $\alpha$.
Consequently we have $G\cong G^\alpha$ for every $\alpha\in{\rm Ant}(G)$, so $G$
is neighborhood reconstructible and a cancellation graph by
Corollaries~\ref{Cor:NbhdAnti} and~\ref{COR:can}.
\hfill \rule{7pt}{7pt}

\medskip
As an example of Proposition~\ref{PROP:biprevinv}, the graph in
Figure~\ref{FigReconstruct} is
neighborhood reconstructible and a cancellation graph.

Corollary~\ref{COR:can2} and Proposition~\ref{PROP:biprevinv} use the
absence of certain kinds of involutions to draw conclusions about whether
a graph is a cancellation graph (neighborhood reconstructible).
An interesting problem would be to find a way to extend the sufficient condition
of Corollary~\ref{COR:can2} to some kind of characterization, as in
Proposition~\ref{PROP:biprevinv}. We leave this as an open problem. 


\end{document}